\documentclass[a4paper,12pt]{article}
\usepackage{amsmath, amsfonts, amssymb}

\makeatletter
\newbox\bk@bxb
\newbox\bk@bxa
\newif\if@bkcont
\newcount\bk@lcnt

\def\breakboxskip{2pt}
\def\breakboxparindent{1.8em}

\def\breakbox{\vskip\breakboxskip\relax
\setbox\bk@bxb\vbox\bgroup
\advance\linewidth -2\fboxrule
\hsize\linewidth\@parboxrestore
\parindent\breakboxparindent\relax}

\def\bk@split{%
\@tempdimb\ht\bk@bxb 
\advance\@tempdimb\dp\bk@bxb
\setbox\bk@bxa\vsplit\bk@bxb to\z@ 
\setbox\bk@bxa\vbox{\unvbox\bk@bxa}
\setbox\@tempboxa\vbox{\copy\bk@bxa\copy\bk@bxb}
\advance\@tempdimb-\ht\@tempboxa
\advance\@tempdimb-\dp\@tempboxa}

\def\bk@addfsepht{%
\setbox\bk@bxa\vbox{\vskip\fboxsep\box\bk@bxa}}

\def\bk@addskipht{%
\setbox\bk@bxa\vbox{\vskip\@tempdimb\box\bk@bxa}}

\def\bk@addfsepdp{%
\@tempdima\dp\bk@bxa
\advance\@tempdima\fboxsep
\dp\bk@bxa\@tempdima}

\def\bk@addskipdp{%
\@tempdima\dp\bk@bxa
\advance\@tempdima\@tempdimb
\dp\bk@bxa\@tempdima}

\def\bk@line{%
\hbox to \linewidth{%
\hskip-2\fboxsep\vrule \@width\fboxrule\hskip.5\fboxsep\vrule \@width\fboxrule\hskip1.5\fboxsep
\box\bk@bxa\hfil
}}%

\def\endbreakbox{\egroup
\ifhmode\par\fi{\noindent\bk@lcnt\@ne
\@bkconttrue\baselineskip\z@\lineskiplimit\z@
\lineskip\z@\vfuzz\maxdimen
\bk@split\bk@addfsepht\bk@addskipdp
\ifvoid\bk@bxb 
\def\bk@fstln{\bk@addfsepdp
\hskip-\parindent\vbox{\llap{\raisebox{-2ex}{\rule{1.5\fboxsep}{\fboxrule}\hskip.5\fboxsep}}\bk@line\llap{\rule{1.5\fboxsep}{\fboxrule}\hskip.5\fboxsep}}}

\else 
\def\bk@fstln{\vbox{\llap{\raisebox{-2ex}{\rule{1.5\fboxsep}{\fboxrule}\hskip.5\fboxsep}}\bk@line}\hfil%
\advance\bk@lcnt\@ne
\loop
\bk@split\bk@addskipdp\leavevmode
\ifvoid\bk@bxb 
\@bkcontfalse\bk@addfsepdp
\vtop{\bk@line\noindent\hskip-2\fboxsep{\rule{1.5\fboxsep}{\fboxrule}}}%

\else 
\bk@line
\fi
\hfil\advance\bk@lcnt\@ne
\if@bkcont\repeat}%
\fi
\leavevmode\bk@fstln\par}\vskip\breakboxskip\relax}

\def\smp{\smallskip\par}
\def\un{{\bf 1}}

\def\pf{\noindent{\bf Proof~:}\ }
\def\findemo{~\leaders\hbox to 1em{\hss\  \hss}\hfill~\raisebox{.5ex}{\framebox[1ex]{}}\smp}

\def\mpn{\medskip\par\noindent}
\def\smpn{\smallskip\par\noindent}
\def\normal{\mathop{\trianglelefteq}}

\def\smp{\smallskip\par}
\def\smpn{\smallskip\par\noindent}

\def\mpoint{\;\;.}
\def\mvirg{\;\;,}

\def\Res{{\rm Res}}
\def\Ind{{\rm Ind}}
\def\Inf{{\rm Inf}}
\def\Def{{\rm Def}}
\def\Iso{{\rm Iso}}

\def\Inf{{\rm Inf}}

\def\Ker{{\rm Ker}}

\def\Z{\mathbb{Z}}

\newcommand{\romain}[1]{\uppercase\expandafter{\romannumeral #1}}

\newcommand{\flh}[2]{\mathop{\hbox to 12mm{\rightarrowfill}}_{\displaystyle #2}^{\displaystyle #1}\limits}
\newcommand{\sflh}[2]{\mathop{\hbox to 12mm{\rightarrowfill}}_{\scriptstyle #2}^{\scriptstyle #1}\limits}

\newcommand{\sumb}[2]{\mathop{\sum}_{{\scriptstyle #1}\atop {\scriptstyle #2}}}

\newcommand{\carre}[8]{\begin{array}{ccc}
#1&\mathop{\hbox to 12mm{\rightarrowfill}}^{\displaystyle{#2}}\limits&#3\\
\llap{$\displaystyle{#4}$}\left\downarrow\vbox to 6mm{}\right. & & \left\downarrow\vbox to 6mm{}\right.\rlap{$\displaystyle{#5}$}\\
#6&\mathop{\hbox to 12mm{\rightarrowfill}}_{\displaystyle #7}\limits&#8\\
\end{array}}

\newcommand{\carrem}[8]{\begin{array}{ccc}
#1&\mathop{\hbox to 12mm{\rightarrowfill}}^{\displaystyle #2}\limits&#3\\
\llap{$\displaystyle #4$}\left\uparrow\vbox to 6mm{}\right. & & \left\uparrow\vbox to 6mm{}\right.\rlap{$\displaystyle #5$}\\
#6&\mathop{\hbox to 12mm{\rightarrowfill}}_{\displaystyle #7}\limits&#8\\
\end{array}}

\newenvironment{enonce}[1]{\pagebreak[2]\refstepcounter{subsection}\refstepcounter{prop}\smpn{{\bf \thesection.\arabic{prop}.\ \ #1~:}}\begin{it} }{\end{it}\smp}
\newenvironment{enonce*}[1]{\pagebreak[2]\smpn{#1~:}\begin{it} }{\end{it}\smp}
\newcommand{\result}[1]{\begin{enonce}{#1}}
\def\fresult{\end{enonce}}

\newenvironment{mth}[1]{\begin{breakbox}\begin{enonce}{#1}}{\end{enonce}\end{breakbox}}
\newenvironment{mth*}[1]{\begin{breakbox}\begin{enonce*}{#1}}{\end{enonce*}\end{breakbox}}
\newenvironment{rem}[1]{\refstepcounter{subsection}\refstepcounter{prop} \mpn{{\bf \thesection.\arabic{prop}.}\ \ \bf#1\ :}}{\smp}

\def\dom{\backslash}
\makeatletter
\renewenvironment{enumerate}{\ifnum \@enumdepth >3 \@toodeep\else
      \advance\@enumdepth \@ne
      \edef\@enumctr{enum\romannumeral\the\@enumdepth}\list
      {\csname label\@enumctr\endcsname}{\setlength{\topsep}{1ex}\setlength{\itemsep}{0pt}\usecounter
        {\@enumctr}\def\makelabel##1{\hss\llap{##1}}}\fi}{\endlist}

\makeatother
\makeatletter
\def\@sect#1#2#3#4#5#6[#7]#8{\ifnum #2>\c@secnumdepth
    \let\@svsec\@empty\else
    \refstepcounter{#1}\edef\@svsec{\csname the#1\endcsname .\hskip .5em}\fi
    \@tempskipa #5\relax
     \ifdim \@tempskipa>\z@
       \begingroup #6\relax
         \@hangfrom{\hskip #3\relax\@svsec}{\interlinepenalty \@M #8\par}%
       \endgroup
      \csname #1mark\endcsname{#7}\addcontentsline
        {toc}{#1}{\ifnum #2>\c@secnumdepth \else
                     \protect\numberline{\csname the#1\endcsname}\fi
                   #7}\else
       \def\@svsechd{#6\hskip #3\relax  
                  \@svsec #8\csname #1mark\endcsname
                     {#7}\addcontentsline
                          {toc}{#1}{\ifnum #2>\c@secnumdepth \else
                            \protect\numberline{\csname the#1\endcsname}\fi
                      #7}}\fi
    \@xsect{#5}}
\def\section{\@startsection {section}{1}{\z@}{-3.5ex plus-1ex minus
    -.2ex}{2.3ex plus.2ex}{\reset@font\Large\bf}}  

\makeatother

\def\mar[#1]{\ar@{-}[#1]|-{\object@{<}}}
\def\marb[#1]{\ar@{-}[#1]|{\object+{  }}}

\def\KB{\mathbb{K}B}
\def\K{\mathbb{K}}
\begin{document}

\centerline{\Large\bf A conjecture on $B$-groups}\vspace{.5cm}\par
\centerline{\bf Serge Bouc}\vspace{1cm}\par
{\footnotesize {\bf Abstract :}  In this note, I propose the following conjecture~: a finite group $G$ is nilpotent if and only if its largest quotient $B$-group $\beta(G)$ is nilpotent. I give a proof of this conjecture under the additional assumption that $G$ be solvable. I also show that this conjecture is equivalent to the following~: the kernel of restrictions to nilpotent subgroups is a biset-subfunctor of the Burnside functor.\vspace{2ex}\par}
{\footnotesize {\bf AMS Subject classification :} 18B99, 19A22, 20J15.\vspace{2ex}}\par

\section{Introduction}
In the study of the lattice of biset-subfunctors of the Burnside functor $\KB$ over a field $\mathbb{K}$ of characteristic 0 (cf. Section~7.2 of \cite{doublact2}, or Chapter~5 of~\cite{bisetfunctors}), a special class of finite groups, called {\em $B$-groups}, plays an important role. It was shown in particular in~\cite{doublact2} Proposition~9 (see also~\cite{bisetfunctors} Theorem 5.4.11) that any finite group $G$ admits a largest quotient in this class, well-defined up to isomorphism, and denoted by $\beta(G)$. A few properties of $B$-groups were proved in~\cite{doublact2}, some of which will be recalled in this paper, but almost no progress was made since, until the following theorem proved recently by M\'elanie Baumann (\cite{baumann-jofa})~: when $p$ is a prime number, recall that a finite group $G$ is called {\em cyclic modulo $p$} (or {\em $p$-hypo-elementary}) if the group $G/O_p(G)$ is cyclic.
\begin{mth}{Theorem} {\rm [M. Baumann]} \label{baumann}Let $p$ be a prime number, and $G$ be a finite group. Then $G$ is cyclic modulo $p$ if and only if $\beta(G)$ is cyclic modulo~$p$.
\end{mth}
In this note I propose the following similar looking conjecture~:
\begin{mth*}{{\bf Conjecture A}} Let $G$ be a finite group. Then $\beta(G)$ is nilpotent if and only if $G$ is nilpotent.
\end{mth*}
\begin{rem}{Remark} \label{nilpotent B-group}It was shown in~\cite{doublact} (Proposition 14) that the nilpotent $B$-groups are the groups of the form $C_n\times C_n$, where $C_n$ is a cyclic group of {\em square free} order~$n$.
\end{rem}
After recalling the basic definitions and properties of $B$-groups, I will give a proof of Conjecture A under the additional assumption that $G$ be {\em solvable}.
\section{$B$-groups}
Let $\mathbb{K}$ be a field of characteristic 0. Let $G$ be a finite group, let $s_G$ denote the set of subgroups of $G$, and let $[s_G]$ be a set of representatives of $G$-conjugacy classes on $s_G$. \par
Denote by $\KB(G)$ the Burnside algebra of $G$ over $\mathbb{K}$. It is a split semisimple commutative $\mathbb{K}$-algebra, with two natural $\K$-bases~: the first one consists of the isomorphism classes of transitive $G$-sets, i.e. the set $\{[G/H]|H\in[s_G]\}$. The second one consists of the primitive idempotents of $\KB(G)$, i.e. the set $\{e_H^G|\ H\in[s_G]\}$. The transition matrix from the first basis to the second one has been described explicitly by Gluck~(\cite{gluck}) and Yoshida~(\cite{yoshidaidemp}), as follows
$$e_H^G=\frac{1}{N_G(H)|}\sum_{X\leq H}|X|\mu(X,H)\,[G/X]\mvirg$$
where $\mu$ is the M\"obius function of the poset of subgroups of $G$.\par
The correspondence $G\mapsto \KB(G)$ is a biset functor~: when $G$ and $H$ are finite groups, and $U$ is a finite $(H,G)$-biset, the functor $S\mapsto U\times_GS$ from the category of finite $G$-sets to the category of finite $H$-sets induces a map $\KB(U):\KB(G)\to \KB(H)$, which is well behaved with respect to disjoint union and composition of bisets.\par
This involves in a single formalism the usual operations of restriction, induction, inflation, and transport by isomorphism between the corresponding Burnside groups. It also involves the less usual operation of {\em deflation}~: when $N\normal G$, the deflation homomorphism $\Def_{G/N}^G:B(G)\to B(G/N)$ corresponds to the $(G/N,G)$-biset $G/N$, and it induced by the functor $S\mapsto N\dom S$ from $G$-sets to $(G/N)$-sets.\par
These elementary operations can be expressed explicitly in each of the above bases. In particular (\cite{bisetfunctors} Theorem 5.2.4)~:
\begin{mth}{Theorem}  Let $G$ be a finite group.
\begin{enumerate}
\item If $H$ is a proper subgroup of $G$, then $\Res_H^Ge_G^G=0$.
\item When $N\normal G$, set
$$m_{G,N}=\frac{1}{|G|}\sumb{X\leq G}{\rule{0ex}{1.3ex}XN=G}|X|\mu(X,G)\mpoint$$
Then $\Def_{G/N}^Ge_G^G=m_{G,N}e_{G/N}^{G/N}$.
\end{enumerate}
\end{mth}
This leads to the notion of $B$-group~: the group $G$ is a $B$-group if any proper deflation of $e_G^G$ is equal to 0. In other words~:
\begin{mth}{Definition} The finite group $G$ is called a $B$-group if $m_{G,N}=0$ for any non-trivial normal subgroup $N$ of $G$.
\end{mth}
\begin{mth}{Notation} When $G$ is a finite group, and $N\normal G$ is maximal such that $m_{G,N}\neq 0$, set $\beta(G)=G/N$.
\end{mth}
There may be several normal subgroups $N$ with the required properties, but the group $G/N$ does not depend on the choice of $N$, up to isomorphism. More precisely (\cite{bisetfunctors} Theorem 5.4.11)~:
\begin{mth}{Theorem} Let $G$ be a finite group.
\begin{enumerate}
\item The group $\beta(G)$ is a $B$-group.
\item If a $B$-group $H$ is isomorphic to a quotient of $G$, then $H$ is isomorphic to a quotient of $\beta(G)$.
\item Let $M\normal G$. The following conditions are equivalent~:
\begin{enumerate}
\item $m_{G,M}\neq 0$.
\item The group $\beta(G)$ is isomorphic to a quotient of $G/M$.
\item $\beta(G)\cong\beta(G/M)$.
\end{enumerate}
\end{enumerate}
\end{mth}
\begin{mth}{Proposition} Let $G$ be a finite group.
\begin{enumerate}
\item The group $G$ is a $B$-group if and only if $m_{G,N}=0$ for any minimal (non-trivial) normal subgroup of $G$.
\item Let $N$ be a minimal (non-trivial) normal subgroup of $G$. If $N$ is abelian, then
$$m_{G,N}=1-\frac{|K_G(N)|}{|N|}\mvirg$$
where $K_G(N)$ is the set of complements of $N$ in $G$.
\end{enumerate}
\end{mth}
\pf Assertion~1 follows from the transitivity of deflations. Assertion~2 is Proposition~5.6.4 of~\cite{bisetfunctors}.\findemo
\section{Proof of Conjecture A in the solvable case}
\begin{mth}{Theorem} \label{solvable case}Let $G$ be a solvable finite group. Then $\beta(G)$ is nilpotent if and only if $G$ is nilpotent.
\end{mth}
\pf If $G$ is nilpotent, then $\beta(G)$ is nilpotent, for it is a quotient of~$G$. The converse follows from an induction argument on the order of $G$~: assume that if $G'$ is a finite solvable group of order $|G'|<|G|$, and if $\beta(G')$ is nilpotent, then $G'$ is nipotent. Assume that $\beta(G)$ is nilpotent, and let $N$ be a non-trivial normal subgroup of $G$. Since $\beta(G/N)$ is a quotient of $\beta(G)$, it is nilpotent. Hence $G/N$ is nilpotent. In particular, if $Z(G)\neq \un$, then $G/Z(G)$ is nilpotent, hence $G$ is nilpotent. So we can assume that $Z(G)=\un$.\par
Now suppose that $M$ and $N$ are non trivial normal subgroups of $G$, such that $M\cap N=\un$. Then $G$ is nilpotent~: indeed, the group $G/(M\cap N)$, isomorphic to $G$, maps injectively into $(G/M)\times (G/N)$, which is nilpotent. \par
It follows that we can assume that $G$ has a unique (non trivial) minimal normal subgroup $N$. Since $G$ is solvable, the group $N$ is elementary abelian, isomorphic to $(C_p)^k$, for some prime number $p$ and some integer $k\geq 1$. Let $Q$ be a Sylow $p$-subgroup of $G$. Then $Q\geq N$. Since the group $G/N$ is nilpotent, and since $Q/N$ is a Sylow $p$-subgroup of $G/N$, it follows that $Q/N\normal G/N$, i.e. $Q\normal G$. \par
Now $N$  is a non trivial normal subgroup of $Q$, thus $N\cap Z(Q)\neq\un$. But $N\cap Z(Q)$ is a normal subgroup of $G$, and by minimality of $N$, it follows that $N\leq Z(Q)$.\par
The group $G$ splits as a semidirect product $G=Q\rtimes H$, where $H$ is a (nilpotent) $p'$-subgroup of $G$. The group $H$ acts on the $p$-group $Q$, thus $Q=C_Q(H)[H,Q]$ (by \cite{gorenstein} Theorem 3.5 Chapter 5). \par
Since $G/N$ is nilpotent, it follows that $G/N\cong (Q/N)\times H$. It follows that $[H,Q]\leq N$. Thus $Q=C_Q(H)N$. Now $N\cap C_Q(H)$ is centralized by $H$, and by $Q$, since $N\leq Z(Q)$. Thus $N\cap C_Q(H)\leq Z(G)=\un$, and it follows that $Q=N\times C_Q(H)$. Then $C_Q(H)$ is normalized by $Q$, and centralized by~$H$. Thus $C_Q(H)\normal G$, and as $N\cap C_Q(H)=\un$, it follows that $C_Q(H)=\un$, thus $N=Q$.\par
But now $G=N\rtimes H$, where $N\cong (C_p)^k$, and $H$ is a $p'$-group. Since $N$ is minimal normal in $G$, it follows that $H$ acts irreducibly on $N$, and that $H$ is a maximal subgroup of $G$. Since $H$ is not normal in $G$ (as $N$ is the only minimal normal subgroup of $G$, and $N\nleq H$), it follows that $N_G(H)=H$.
Finally, since $H$ is a $p'$-group, all the complements of $N$ in $G$ are conjugate, hence $|K_G(N)|=|G:N_G(H)|=|N|$. Thus $m_{G,N}=1-\frac{|K_G(N)|}{|N|}=0$. Hence $G$ is a $B$-group, thus $G\cong \beta(G)$ is nilpotent.\findemo
\begin{rem}{Remark} Actually, in the situation of the end of the proof, the group $G$ is trivial~: indeed, it is a nilpotent $B$-group, hence isomorphic to $C_n\times C_n$, where $n$ is a square free integer. As $G$ has a unique minimal normal subgroup by assumption, the only possibility is $n=1$.
\end{rem}
\section{\bf Comments}
The following conjecture doesn't mention $B$-groups~:
\begin{mth*}{{\bf Conjecture B}} For any group $G$, let $\nu_G$ denote the restriction map
$$\nu_G=\prod_{H\in\mathcal{N}(G)}\Res_H^G :B(G)\to\prod_{H\in \mathcal{N}(G)}B(H)\mvirg$$
where $\mathcal{N}(G)$ is the set of nilpotent subgroups of $G$.\par
Then the correspondence $G\mapsto L(G)=\Ker\,\nu_G$ is a biset subfunctor of~$B$.
\end{mth*}
Still~:
\begin{mth}{Theorem} \label{B eqv A}Conjecture B is equivalent to Conjecture A.
\end{mth}
\pf Since $B(G)$ is a free $\Z$-module, it maps injectively in $\KB(G)$. Let $u\in B(G)$. Then $u$ can be written
$$u=\sum_{H\in[s_G]} |u^H|e_H^G$$
in $\KB(G)$. Thus $u\in L(G)$ if and only if $|u^H|=0$ for any $H\in\mathcal{N}(G)$.\par
Suppose that Conjecture A holds. Proving Conjecture B amounts to proving that $L$ is invariant under the elementary biset operations of induction, restriction, inflation, deflation, and transport by isomorphism. \par
The latter case is clear~: if $\varphi: G\to G'$ is a group isomorphism, then $\Iso(\varphi)\big(L(G)\big)\leq L(G')$.\par
Now let $G$ be a group, and let $K$ be a subgroup of $G$. As nilpotent subgroups of $K$ are nilpotent subgroups of $G$, the transitivity of restrictions implies that $\Res_K^GL(G)\subseteq L(K)$. Conversely, if $u\in L(K)$ and $H\in \mathcal{N}(G)$, then by the Mackey formula
$$\Res_{H}^G\Ind_{K}^Gu=\sum_{g\in[H\dom G/K]}\Ind_{H\cap{^gK}}^Hc_g\Res_{H^g\cap K}^Ku=0\mvirg$$
(where $c_g$ denote conjugation by $G$), since $H^g\cap K\in\mathcal{N}(K)$. It follows that $\Ind_K^GL(K)\subseteq L(G)$.\par
Suppose now that $N\normal G$, and that $u\in L(G/N)$. Then for any $H\in\mathcal{N}(G)$
$$\Res_H^G\Inf_{G/N}^Gu=\Inf_{H/H\cap N}^H\Iso_{HN/N}^{H/H\cap N}\Res_{HN/N}^{G/N}u=0\mvirg$$
since $HN/N\in\mathcal{N}(G/N)$. Hence $\Inf_{G/N}^GL(G/N)\subseteq L(G)$.\par
Finally, in the same situation, let $v\in L(G)$, and $K/N\in\mathcal{N}(G/N)$. Then
$$\Res_{K/N}^{G/N}\Def_{G/N}^Gv=\Def_{K/N}^K\Res_K^Gv\mpoint$$
Moreover
$$\Res_K^Gv=\sum_{X\in[s_K]}|v^X|e_X^K\mvirg$$
and $|v^X|=0$ for $X\in\mathcal{N}(X)$. Since moreover
$$e_X^K=\frac{1}{|N_K(X):X|}\Ind_X^Ke_X^X\mvirg$$
 it follows that
\begin{eqnarray*}
\Res_{K/N}^{G/N}\Def_{G/N}^Gv&=&\sum_{X\in[s_K]}\frac{|v^X|}{|N_K(X):X|}\Def_{K/N}^K\Ind_X^Ke_X^X\\
&=&\sum_{X\in[s_K]}\frac{|v^X|}{|N_K(X):X|}\Ind_{XN/N}^{K/N}\Iso_{X/X\cap N}^{XN/N}\Def_{X/X\cap N}^Xe_X^X\\
&=&\sum_{X\in[s_K]}\frac{|v^X|}{|N_K(X):X|}m_{X,X\cap N}\Ind_{XN/N}^{K/N}\Iso_{X/X\cap N}^{XN/N}e_{X/X\cap N}^{X/X\cap N}\\
&=&\sum_{X\in[s_K]}\frac{|v^X|}{|N_K(X):X|}m_{X,X\cap N}\Ind_{XN/N}^{K/N}e_{XN/N}^{XN/N}\mpoint
\end{eqnarray*}
If $X\leq K$ is such that $m_{X,X\cap N}\neq 0$, then $\beta(X)\cong\beta(X/X\cap N)\cong\beta(XN/N)$. The group $XN/N$ is a subgroup of $K/N$, hence it is nilpotent, hence $\beta(X)\cong\beta(XN/N)$ is nilpotent. If Conjecture~A is true, then $X$ is nilpotent, hence $|v^X|=0$. It follows that $\Res_{K/N}^{G/N}\Def_{G/N}^Gv=0$, hence $\Def_{G/N}^GL(G)\leq L(G/N)$.\par
Observe that one can still conclude that $\Res_{K/N}^{G/N}\Def_{G/N}^Gv=0$ without assuming Conjecture A, in the case where $N$ is solvable~: indeed in this case, the group $K$ is solvable (as $N$ is solvable and $K/N$ is nilpotent), so $X$ is solvable, and one can conclude by Theorem~\ref{solvable case}.\par
Conversely, assume that Conjecture B holds, and let $G$ be a finite group. Then $|G|e_G^G$ is an element of $B(G)$, whose restrictions to all proper subgroups of $G$ are 0. If $G$ is not nilpotent, then $|G|e_G^G\in L(G)$. Hence for any normal subgroup $N$ of $G$, the element 
$$\Def_{G/N}^G|G|e_G^G=m_{G,N}|G|e_{G/N}^{G/N}$$
is in $L(G/N)$. If $G/N\cong \beta(G)$, then $m_{G,N}\neq 0$, hence $|G|e_{G/N}^{G/N}\in L(G/N)$. Since it is a non-zero element, it follows that $G/N\notin\mathcal{N}(G/N)$, i.e. that $G/N\cong \beta(G)$ is not nilpotent. Hence Conjecture A holds.\findemo
\begin{rem}{Remark} The above proof shows that the correspondence $G\mapsto L(G)$ is a biset functor on the full subcategory of the biset category consisting of {\em solvable} groups. Actually, it proves a little more~: the correspondence $G\mapsto L(G)$ is a biset functor on the category of {\em all} finite groups, if we only allow as morphisms those bisets for which {\em left stabilizers are solvable} (or equivalently, if we only allow deflations by solvable normal subgroups).
\end{rem}
\begin{rem}{Remark} Let $G$ be a minimal counterexample to Conjecture~A. Then $G$ is non solvable, and as above $G$ has a unique minimal normal subgroup $N$, non-central in $G$. Thus $N\cong S^k$, where $S$ is a non-abelian simple group, and $C_G(N)=\un$. The group $H=G/N$ is nilpotent, and it is the largest solvable quotient of $G$. In particular $\beta(G)$ is a quotient of $H$, hence of $\beta(H)$. Since $\beta(G)$ cannot be a $p$-group (for otherwise $G$ would be cyclic modulo $p$, by Theorem~\ref{baumann}, hence solvable, hence nilpotent by Theorem~\ref{solvable case}), it follows that $H$ is not $p$-elementary for any prime $p$ (that is, there are at least two different primes $p$ such that the Sylow $p$-subgroups of $H$ are non cyclic)~: indeed, if $P$ is a $p$-group, then $\beta(P)=\un$ if $P$ is cyclic, and $\beta(P)=C_p\times C_p$ otherwise (cf. Remark~\ref{nilpotent B-group}). \par
Since $\Phi(G)$ is nilpotent, it follows that $\Phi(G)=\un$. Let $X$ be a minimal subgroup of $G$ such that $XN=G$. Then $X\cap N\leq \Phi(X)$, thus $m_{X,X\cap N}=1$. Hence $\beta(X)\cong \beta(X/X\cap N)\cong\beta(G)$ is nilpotent. Moreover $X<G$ (as $N\nleq\Phi(G)=\un$), so $X$ is nilpotent.
\end{rem}
\begin{rem}{Remark} In view of Conjecture A and Theorem~\ref{solvable case}, Jacques Th\'evenaz has proposed the following~:
\begin{mth*}{{\bf Conjecture}} Let $G$ be a finite group. Then $\beta(G)$ is solvable if and only if $G$ is solvable.
\end{mth*}
This conjecture implies Conjecture A, by Theorem~\ref{solvable case}. A straightforward modification of the proof of Theorem~\ref{B eqv A} shows that this conjecture is equivalent to saying that the correspondence sending a finite group $G$ to the kernel of the restriction map 
$$\rho_G=\prod_{H\in\mathcal{R}(G)}\Res_H^G :B(G)\to\prod_{H\in \mathcal{R}(G)}B(H)\mvirg$$
where $\mathcal{R}(G)$ is the set of solvable subgroups of $G$, is a biset functor.
\end{rem}

\centerline{\rule{5ex}{.1ex}}
\begin{flushleft}
Serge Bouc - CNRS-LAMFA, Universit\'e de Picardie, 33 rue St Leu, 80039, Amiens Cedex 01 - France. \\
{\tt email : serge.bouc@u-picardie.fr}\\
{\tt web~~ : http://www.lamfa.u-picardie.fr/bouc/}
\end{flushleft}
\end{document}